\documentclass[12pt]{amsart}

\usepackage{amsfonts, amssymb, curves, epic, epsf, amsmath}
\newcommand{\R}{\mathbb{R}}

\newcommand{\Z}{\mathbb{Z}}
\newcommand{\cT}{\mathbb{T}}

\newcommand{\cP}{\mathcal{P}}

\newcommand{\Cech}{{\v Cech}{} }
\newcommand{\Om}{\Omega}
\newcommand{\cU}{\mathcal{U}}
\newcommand{\cV}{\mathcal{V}}

\newcommand{\tsigma}{\tilde{\sigma}}

\newcommand{\larr}{\left( \begin{array}{c}}
\newcommand{\rarr}{\end{array} \right) }

\newcommand{\lsqarr}{\left[ \begin{array}{c}} 
\newcommand{\rsqarr}{\end{array} \right]}

\newcommand{\arrow}{\rightarrow}

\newcommand{\inv}{\varprojlim}
\newcommand{\dir}{\varinjlim}
\newcommand{\seabox}{{\framebox{$\searrow$}}}
\newcommand{\neabox}{{\framebox{$\nearrow$}}}
\newcommand{\swabox}{{\framebox{$\swarrow$}}}
\newcommand{\nwabox}{{\framebox{$\nwarrow$}}}
\newcommand{\neswabox}{{\framebox{$\nearrow\mkern-18mu\swarrow$}}}
\newcommand{\nwseabox}{{\framebox{$\nwarrow\mkern-18mu\searrow$}}}
\newcommand{\ER}{{\hbox{\tiny ER}}}

\begin{document}

\newtheorem{theorem}{Theorem}
\newtheorem{cor}[theorem]{Corollary}
\newtheorem{corollary}[theorem]{Corollary}
\newtheorem{lemma}[theorem]{Lemma} \newtheorem{prop}[theorem]{Proposition}
\newtheorem{example}[theorem]{Example}

\title{Cohomology of substitution tiling spaces}
\author{Marcy Barge, Beverly Diamond, John Hunton and Lorenzo Sadun}

 \maketitle
\begin{center} May 25, 2009 \end{center}
 \begin{abstract}

   Anderson and Putnam showed that the cohomology of a substitution
   tiling space may be computed by collaring tiles to obtain a
   substitution which ``forces its border.''  One can then represent
   the tiling space as an inverse limit of an inflation and
   substitution map on a cellular complex formed from the collared
   tiles; the cohomology of the tiling space is computed as the direct
   limit of the homomorphism induced by inflation and substitution on
   the cohomology of the complex.  In earlier work, Barge and Diamond
   described a modification of the Anderson-Putnam complex on collared
   tiles for one-dimensional substitution tiling spaces that allows
   for easier computation and provides a means of identifying certain
   special features of the tiling space with particular elements of
   the cohomology.  In this paper, we extend this modified
   construction to higher dimensions. We also examine the action of
   the rotation group on cohomology and compute the cohomology of the
   pinwheel tiling space.
\end{abstract} 


\section{Introduction}

The study of aperiodic tilings of Euclidean space began in 1961 with
the work of Wang \cite{Wang}, whose interest centered on decidability
issues.  Schechtman's discovery of quasicrystaline materials in 1984
\cite{quasi} shifted focus to the combinatorial structure of
nonperiodic tilings. Tiling space now arises naturally: by passing to
the space of all tilings locally indistinguishable from a given
tiling, combinatorial properties of an individual tiling are
translated into topological properties of the tiling space. This space
can be studied with the aid of naturally defined dynamical systems.

There are two main approaches to the study of tiling dynamics.
One approach associates a C*-algebra to a group action on the tiling space 
and then studies the K-theory of this algebra. In the topological approach,
cohomology is computed in some more direct manner. The C*-algebra
approach stays closer to the physics, while the cohomology calculations are
generally more straightforward. At any rate, the K-theory and cohomology typically 
coincide and this permits cohomological interpretation of such physically relevant
notions as gap-labeling (\cite{B},\cite{Be}, \cite{KP}, \cite{S}). 

The two most general procedures for constructing nonperiodic tilings are
the cut-and-project method and the substitution method. For the calculation
of cohomology of cut-and-project spaces see \cite{FHK}, \cite{GHK}. We consider only
substitution tilings in this article.

In 1998, Anderson and Putnam \cite{ap} showed how to describe a
substitution tiling space as an inverse limit of branched manifolds.
There are two models, depending on whether the substitution has a
property called ``forcing the border''. If the substitution has that
property, they construct a CW complex
$K$ consisting of one copy of every kind of tile, with certain edge 
identifications. The substitution $\sigma$ maps $K$ to itself and the
inverse limit  $\inv (K,\sigma)$ is homeomorphic to the tiling space
$\Omega_\sigma$. 

If the substitution doesn't force the border, then $\inv (K,\sigma)$
is still well-defined, and there is a natural map $\Omega_\sigma \to
\inv (K,\sigma)$, but this map fails to be injective.  In that case, they
build a complex $K_c$ from multiple copies of each tile type, one copy
for each pattern of nearest-neighbor tiles that can touch the tile in
question. A tile, together with a label indicating the pattern of
neighboring tiles, is called a ``collared tile''. Then $K_c$ is obtained
by identifying edges of collared tiles in a particular way, and
$\Omega_\sigma$ is always homeomorphic to $\inv (K_c, \sigma)$. 
(In \cite{k} Kellendonk indepedently introduced collaring to compute the
top dimensional cohomology of some tiling spaces.)

Anderson-Putnam collaring is based on labeling {\em tiles}. In this paper
we extend a 1-dimensional construction of Barge and Diamond \cite{bd1} to
develop a collaring scheme based on labeling {\em points} by their
neighborhoods to a distance $t$. These labeled points naturally
aggregate into a branched manifold $K_t$, and we show how to construct
tiling spaces (not just substitution tiling spaces) as inverse limits
of the branched manifolds $K_t$.  Moreover, our construction extends to
tiling spaces with continuous rotational symmetry (e.g., the pinwheel
tiling), and to spaces of tilings that do not have finite local
complexity (see \cite{FS}).

There are two important inverse limits.  The first applies to all
tiling spaces, whether arising via substitutions or not. If $t' > t$, then there is a
natural map $f: K_{t'} \to K_t$ that merely forgets collaring
data from farther than a distance $t$. The tiling space $\Omega$ is
always homeomorphic to $\inv (K_t, f)$. 

The second inverse limit applies to substitution tiling spaces.  If
$\sigma$ is an expanding substitution, then there is a constant
$\lambda>1$ such that $\sigma$ maps $K_{t'}$ to all $K_t$ with
$t<\lambda t'$. In particular, we can take $t'=t$ and compute
$\inv(K_t,\sigma)$. For any positive $t$, $\Omega_\sigma$ is
homeomorphic to this tiling space, and the \Cech cohomology
of $\Omega_\sigma$ is computed as $\check H^*(\Omega_\sigma) = \dir
(H^*(K_t),\sigma^*)$.

The problem with this construction is that $\sigma$ does not respect
the cellular structure of $K_t$, so computing $\sigma^*$ can be
difficult.  The solution is to take a cellular map $\tsigma:K_t \to K_t$,
homotopic to $\sigma$, and consider the space $\Xi = \inv(K_t,
\sigma^*)$.  In general, $\Xi$ will not be homeomorphic to $\Omega_\sigma$,
but they will both have (\Cech) cohomology $\dir (H^*(K_t), \tsigma^*)$.

An added benefit is stratification. There are many stratifications of 
$K_t$ (i.e., subsets $S_0 \subset S_1 \subset S_2 \subset \cdots \subset 
K_t$) and we can often pick a $\tsigma$ that maps each stratum $S_k$ to
itself. We then consider the inverse limits $\Xi_k = \inv (S_k,
\tsigma)$ and compute the relative cohomology groups $\check
H^*(\Xi_{k+1},\Xi_k) = 
\dir (H^*(S_{k+1},S_k),\tsigma^*)$. These assist both in computing 
$\check H^*(\Xi)$ and in interpreting the different terms that appear in the 
final answer. 

In the next section we fix terminology. Section 3 describes the new version
of collaring that we employ in subsequent computations and in Section 4
we set up the stratification of cohomology that this collaring permits.
As a first example, the cohomology of the chair tiling space is computed in Section 5.
Our approach gives an efficient ``by hand'' computation in which topological
features of this space are reflected in the stratified cohomology.

In Section 6 we extend consideration to the space obtained by allowing 
the full Euclidean group to act on a tiling. For a translationally finite tiling
(like a chair tiling or a Penrose tiling) three spaces arise: the closure of the
orbit under translations ($\Omega^1$), the closure of the orbit under the full
Euclidean group ($\Omega^{rot}$), and $\Omega^{rot}$ mod rotations ($\Omega^0$).
We establish general relations between the cohomologies of these spaces and,
by way of example, calculate these for the chair and Penrose. The top
dimensional cohomology of $\Omega^{rot}$ for the Penrose tiling has
5-torsion, illustrating a general result linking $n$-torsion to the existence of 
more than one tiling with $n$-fold rotational symmetry.

In the concluding Section 7 we give the first computation of the cohomology of the
pinwheel tiling space.

\section{Notation and Terminology}

Let $\cP$ be a finite set of compact subsets of $\R^d$, each the closure
of its interior; these subsets shall be called {\em prototiles}.  A {\em
tile\/} is a set obtained from a prototile by rigid motion (by which 
we will mean either a translation or an arbitrary Euclidean
motion, depending on context). A {\em patch for
$\cP$\/} is a set of tiles with pairwise disjoint interiors and the {\em
support} of a patch is the union of its tiles. A {\em tiling of $\R^d$
with prototiles $\cP$} is a patch with support $\R^d$.  If $T$ is a tiling
and $A$ is a bounded subset of
$\R^d$, denote by $[A]_T$ the set of all tiles in $T$ that have
nonempty intersection with $A$.  Two patches $P^1$ and $P^2$ are {\em
translationally equivalent\/} if there is a $w \in \R^d$ so that $P^1 +
w = P^2$.  A tiling has {\em translationally finite local complexity\/},
or is {\em 'translationally finite'}, for short, if, for each $r >
0$, the tiling contains only finitely many translational equivalence
classes of patches of diameter less than $r$.  Two patches $P^1$ and
$P^2$ are {\em rigid equivalent\/} if there is a rigid motion taking
$P^1$ to $P^2$, and a tiling has {\em finite local complexity\/} if, for
each $r>0$, the tiling contains only finitely many rigid equivalence
classes of patches of diameter less than $r$.

We use a topology where two tilings are close if they agree on a large
ball around the origin, up to a small motion of each tile. If the
tilings have finite local complexity, a small motion of each tile must
come from a small rigid motion of the entire tiling.  If the tilings
are translationally finite, this rigid motion must be a
translation.

We first consider substitutions on translationally finite tilings and then
extend the definitions to cover other cases. 
A substitution $\sigma$ on $\cP$ is a map $\sigma: \cP \arrow \cP^*$,
where $\cP^*$ is the collection of finite patches from $\cP$, such
that, for $p \in \cP$, the support of $\sigma(p)$ is the rescaled
prototile $L p$, where $L$ is a fixed expanding linear map (that is,
all eigenvalues of $L$ have modulus larger than one).  The
substitution map can be extended to patches, dilating the entire patch
by a factor of $L$ and replacing each dilated tile $L t_i = L(p + w)$
with $\sigma(t_i) := \sigma(p) + L w $.  A tiling $T$ of $\R^d$ is {\em
admissible} for $\sigma$ if, for every finite patch $P$ of $T$, there
is a prototile $p$ and an integer $n$ such that $P$ is equivalent to
a subpatch of $\sigma^n(p)$.  The {\em substitution tiling space} for
$\sigma$ and $\cP$, written $\Om_\sigma$ or simply $\Om$, is then the set
of all admissible tilings; $\Om$ is also called the {\em continuous hull}.
There are two natural dynamical systems on $\Omega$: translation, and the
action of substitution on entire tilings.

A substitution $\sigma$ is {\em primitive} if there exists an integer
$n$ such that for any two prototiles $p_1$ and $p_2$, $\sigma^n(p_1)$
contains a copy of $p_2$.  In the case $\sigma$ is primitive, $\Om$ is
minimal under translation.  In the following, we assume that $\sigma$
is primitive and $\Om$ is {\em non-periodic}, that is, contains no
tilings periodic under translation.  The substitution $\sigma$ is {\em
translationally finite} (resp., has {\em finite local complexity}) if
every tiling in $\Om$ is translationally finite (resp., has finite
local complexity). Substitutions for which $\sigma: \Om \to \Om$ is a
homeomorphism are called {\em recognizable}.  All non-periodic and
translationally finite substitutions are recognizable (\cite{sol},
\cite{mosse}).

Things are slightly more complicated if we wish to allow rotations.
The substitution maps rotated versions of $p$ to rotated versions of
$\sigma(p)$ and maps rotated versions of $p+w$ to rotated versions of
$\sigma(p)+Lw$.  For this to make sense, $L$ must commute with the
rotations being considered.  In 2 dimensions, this means that $L$ must
be a uniform dilation by a scaling factor $\lambda$, followed by a
rotation. There exist theorems about such substitutions being recognizable
\cite{HRS}, but in 3 or more
dimensions the hypotheses are complicated.

If $\sigma$ is primitive then the space $\Om$ is the closure of the orbit
under translation of any one tiling whose patches are all admissible. If
$\sigma$ is translationally finite, this is the same as the set of
admissible tilings.  If $\sigma$ is not translationally finite, then $\Om$
contains tilings whose patches may not all be admissible, but which can be 
approximated arbitrarily closely by admissible patches. For instance, 
in the pinwheel tiling the admissible patches contain tiles pointing in 
a countable and dense set of directions, but the space $\Om$ allows tiles to 
point in an uncountable continuum of directions. 

Let $X_0, X_1, X_2, \ldots$ be a sequence of compact metric spaces with
continuous {\em bonding maps} $f_n: X_n \to X_{n-1}$. 
The {\em inverse limit\/} $\inv (X,f)$ is the space defined by
\begin{equation}
\{(x_0,x_1,\ldots) \in \prod X_n: f(x_n) = x_{n-1}\;{\rm for}\;
  n=1,2,\ldots\}
\end{equation} 
with the topology inherited from the product topology on
$\prod X_n$. 
The spaces $X_n$ are called {\em approximants} to $\inv (X,f)$. 

For any $x\in \R^d$ and any positive $t$, let $B_t(x)$ denote 
the open ball of radius $t$ around $x$.

\section{The modified complex} \label{complex}

Suppose that $T$ is a translationally finite tiling and pick a positive
real number $t$. We say that two points $x,y \in \R^d$ are {\em equivalent
to distance $t$} if
$[B_t(x)]_T = [B_t(y)]_T + x-y$, and write $x \sim_t y$. 
The space $K_t$ is the quotient of $\R^d$ by this
equivalence relation, with the quotient topology. If $\Omega$ comes from
a primitive substitution, then this quotient space 
will be the same for all choices of $T \in \Omega$ (since all tilings have
the same sets of patches of size $t$), and $K_t$ can be viewed as an 
approximant to the tiling space $\Om$, not just to the tiling $T$.  

To define $K_t$ for tilings that are not necessarily translationally
finite we view $\sim_t$ as an equivalence
relation on the trivial $\R^d$ bundle $E \to \Om$, where the
fibre over $T\in \Om$ may be considered as a copy of $T$ itself.
If $x\in T$ and $y\in T'$, we say that $(x,T) \sim_t (y,T')$ if
$[B_t(x)]_T = [B_t(y)]_{T'} + x-y$.  We then define $K_t$
to be the quotient of $E$ by $\sim_t$.  There is also a natural
projection $\pi: \Om \to K_t$ sending $T$ to the equivalence class of
$0 \in T$.

The space $K_t$ can be viewed as the set of all possible instructions for
tiling a region $B_t(0)$ (and therefore $[B_t(0)]$). 

If $t_2 \ge t_1$, then $x \sim_{t_2} y$ implies $x \sim_{t_1} y$, so there
is a ``forgetful'' map $f: K_{t_2} \to K_{t_1}$. Similarly, if $t$ is
sufficiently large there is a natural quotient map from $K_t$ to the
collared tile complex of Anderson and Putnam \cite{ap}, a result that
extends to maps from $K_t$ for $t$ large to the complex of $n$-fold
collared tiles. We obtain a continuous 
generalization of G\"ahler's construction \cite{gaehler}, applicable to
all tiling spaces.

\begin{theorem}
If $t_0\le t_1\le t_2 \le \ldots$ is an infinite sequence of radii with
$\lim t_n = \infty$, then $\Omega$ is homeomorphic to the inverse limit
$\inv (K_t, f)$
\end{theorem}

\begin{proof}
A point in the inverse limit  $\inv (K_t, f)$
is a set of consistent  instructions for tiling the plane out to all
radii $t_n$, i.e., a set of instructions for tiling the entire plane.
This gives a bijection between
$\inv (K_t, f)$ and $\Om$, and this bijection is easily seen to be a
homeomorphism.  
\end{proof} 

If $L$ is an expansive linear map, then there is a constant $\lambda>1$
and a norm on $\R^d$ such that, for all $t>0$, $B_{\lambda t}(0) \subset
L(B_t(0))$. If $L$ is not diagonalizable, or has eigenvectors that are 
not orthogonal, then this norm may not be the usual Euclidean norm. However,
we can always pick an inner product that is adapted to the geometry of $L$,
or  use linear norms that do not come from an inner product,
such as $L^p$ norms.

Let $\sigma$ be a substitution with expansion $L$. 
If $x \in K_t$, then all tilings in $\pi^{-1}(x)$ agree on 
$B_t(0)$, so all tilings in $\sigma(\pi^{-1}(x))$ agree on $B_{\lambda t}(0)
\subset L(B_t(0))$, so $\sigma$ induces a map $K_t \to K_{\lambda t}$. 
Coupled with the forgetful map, we get maps $\sigma: K_{t_2} \to K_{t_1}$
whenever $t_2 \ge  t_1/\lambda$.

\begin{theorem}\label{limit}
  Let $\sigma$ be a recognizable 
substitution with expansion $L$.  If
  $t_0, t_1, t_2, \ldots$ is a sequence of positive numbers with each
  $t_n \ge t_{n-1}/\lambda$ and with $\lambda^n t_n \to \infty$, then
  the inverse limit  $\inv (K_t,\sigma)$ is homeomorphic to
  $\Om_\sigma$.  
\end{theorem}

\begin{proof} Let $(x_0, x_1, \ldots) \in \inv (K_t, \sigma)$. Each $x_n$
defines a tiling out to distance $t_n$. Applying the substitution $n$ times
gives a tiling on the region $L^n(B_{t_n}(0))$ that, when restricted to 
$L^{n-1}(B_{t_{n-1}}(0))$, agrees with the tiling defined by $x_{n-1}$. 
Since $B_{\lambda^n t_n}(0) \subset L^n(B_{t_n}(0))$ and
$\lambda^n t_n \to \infty$, the point $(x_0, x_1, \ldots)$ defines a 
set of nested patches that exhaust the entire plane. In other words, they
define a tiling. This gives a continuous map 
$\inv (K_t, \sigma) \to \Om_\sigma$. 
The inverse map sends a tiling $T$ to $(x_0, x_1, \ldots)$, where $x_n$ is
the $\sim_{t_n}$ equivalence class of the origin in $\sigma^{-n}(T)$. 
\end{proof}

\begin{corollary}
Let $t>0$, and let $\sigma$ be a recognizable 
substitution with expansion $L$. Then
$\Om_\sigma$ is homeomorphic to $\inv (K_t, \sigma)$, where each approximant
is the same and each map is the same. 
\end{corollary}

The space $K_t$ is easy to visualize if $t$ is much smaller than the
diameter of any tile. Two points, each farther than $t$ from the
boundary of the tiles they sit in, are identified if they sit in
corresponding places in the same tile type.  Points near edges are
identified if they sit in corresponding places in the same tile type,
and the tiles across the nearby edges are of the same type. That is,
these points ``know'' about their neighbor across the edge.  Likewise,
points near vertices ``know'' all the tiles that meet that vertex.

Figure 1 shows the approximants $K_t$ for two 1-dimensional tilings, 
the Thue-Morse substitution $a \to ab, b \to ba$ and the period-doubling
substitution $a \to bb, b \to ba$. In both cases we take the tiles to have
length 1 and take $t < 1/2$. In each case there are intervals $e_a$
and $e_b$ of length
$1-2t$ describing the interiors of the $a$ and $b$ tiles. These are the 
images of tilings that have no vertices in $B_t(0)$.  The tilings with a 
vertex in $B_t(0)$ yield intervals of 
length $2t$, one for each possible transition between tiles. The Thue-Morse
complex has four such ``vertex cells'' $v_{aa}, v_{ab}, v_{ba}, v_{bb}$, 
since the transitions $aa$, $ab$, $ba$, 
and $bb$ are all possible. The period-doubling complex only has three 
such intervals, $v_{ab}$, $v_{ba}$, and $v_{bb}$, 
since consecutive $a$ tiles never occur. 

\begin{figure}
\centerline{\epsfysize=1.8truein\epsfbox{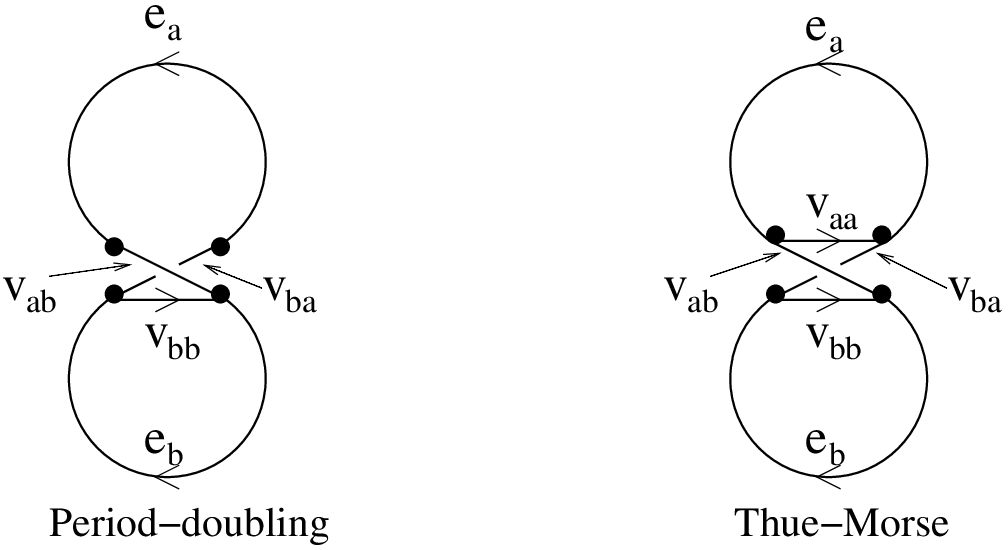}}
\caption{Two approximants}
\label{fig1}
\end{figure}

\section{Homotopy, stratification, and cohomology}

The following theorem about inverse limits in the category of topological
spaces is standard:

\begin{theorem}\label{cch} The \Cech cohomology of an inverse limit 
$\inv (K,f)$ is the direct limit of the \Cech cohomology of the
approximants
$K_n$ under the pullback map $f^*$. If each approximant is a CW complex, 
then this is isomorphic to the direct limit of the 
singular or cellular cohomology of the approximants. 
\end{theorem}

To compute $\check H^*(\Om)$, we need only compute the cellular cohomology
of $K_t$, compute the action of $\sigma^*$ on this cellular cohomology, 
and take the direct limit. Unfortunately, $\sigma$ is not a cellular map. 
In the Thue-Morse and period-doubling examples, the obvious 1-cells in $K_t$ 
are intervals of length
$2t$ and $1-2t$. However, $\sigma$ doubles length, and so sends a cell
of length $2t$ to an interval of length $4t$, that is to a cell of length
$2t$ plus parts of two other cells.  

One solution is to find another map $\tsigma: K_t \to K_t$, that is
cellular and homotopic to $\sigma$. For instance, we can follow the
substitution by a flow that pushes points towards the nearest vertex, so
that points that are $2t$ away from a vertex flow to points just $t$ away
from a vertex.  Then $\check H^*(\Om) = \dir (H^*(K_t), \sigma^*)
= \dir (H^*(K_t),\tsigma^*) = \check H^*(\Xi)$, where $\Xi = \inv (K_t,
\tsigma)$. 

Another, and essentially equivalent, solution when the expansive map $L$ is just magnification by a constant $\lambda>1$ is to note that for all
$t>s>0$ sufficiently small, $K_t$ and $K_s$ are homotopy equivalent, and
thus share the same \Cech cohomology. A corollary of Theorems
\ref{limit} and \ref{cch} tells us that there is an isomorphism
$$\check H^*(\Om)\cong \dir_n (H^*(K_{t\lambda^{-n}}))$$
and note that the map $K_{t\lambda^{-n}}\to K_{t\lambda^{1-n}}$ {\em
is\/} cellular.

In some simple examples, $H^*(K_t)$ and the action of $\tsigma^*$ can be
computed directly. In more complicated examples, it is useful to stratify
$K_t$ into pieces $S_0 \subset S_1 \subset \cdots \subset K_t$ such that 
each stratum $S_k$ is mapped to itself by $\tsigma$. 
Defining $\Xi_k = \inv(S_k,\tsigma)$, we compute $\check H^*(\Xi_0)$ and
the relative groups $\check H^*(\Xi_k, \Xi_{k-1})$ and use the long exact 
sequences of the pair $(\Xi_k,\Xi_{k-1})$ to recursively compute $\check H^*
(\Xi_k)$ for $k>0$.    

In one-dimensional examples, the obvious stratification is given by
taking $S_0$ to be the union of the vertex cells and $S_1$ as $K_t$. 
Barge and Diamond \cite{bd1} used this stratification via the long exact
sequence of the pair $(S_1, S_0)$ to show that $\check H^1(\Om)$ fits
into the exact sequence
\begin{equation} 0 \to \tilde H^0(\Xi_0) \to \dir (\Z^d, M^T) 
\to \check H^1(\Om) \to
\check H^1(\Xi_0) \to 0,\end{equation}
where $d$ is the number of letters and $M$ is the substitution matrix. 
(Here, and in what follows, $\tilde H$ indicates reduced cohomology.)
For instance, in the period-doubling
space, $S_0$ is contractible, so $\tilde H^0(\Xi_0)$ and $\check H^1(\Xi_0)$
vanish, and hence $\check H^1(\Om)$ is the direct limit of the action of
the matrix $M^T = 
\left ( \begin{smallmatrix}
0&2 \cr 1 & 1 \end{smallmatrix} \right )$, 
namely $\Z[1/2] \oplus Z$. In the Thue-Morse
space, $S_0$ has the topology of a circle, on which $\tsigma$ acts 
by reflection.  We then have 
\begin{equation} 0 \to \dir \left (\Z^2, 
\left ( \begin{smallmatrix} 1&1 \cr 1&1\end{smallmatrix}\right ) \right
  ) \to \check H^1(\Om) \to \Z \to 0.\end{equation} 
Once again, $\check H^1(\Om) =
\Z[1/2] \oplus \Z$, only now the factor of $\Z$ comes from the
topology of $S_0$ rather than from the substitution matrix.

In two dimensions, there is some choice over the stratification that can be used. It is often
useful to take $S_0$ to be the points within $t$ of a vertex, $S_1$ the 
points within $t$ of an edge and $S_2 = K_t$, but other stratifications are
also useful: we shall see examples later in our calculations. For
instance, in a tiling by rectangles we could take
$S_0$ to be the points close to a vertex, $S_1$ the points close to a 
horizontal edge, $S_2$ the points close to any edge, and $S_3=K_t$.
We could also take $S_0$ to be the points close to four tiles (i.e., close
to a 4-way crossing), $S_1$ to be the points close to 3 or more tiles
(i.e., close to crossings and T-intersections), $S_2$ to be the points close
to 2 or more tiles (i.e., close to an edge) and $S_3$ to be everything. 
In three or more dimensions, the possibilities are even more varied.  

The ``right'' stratification depends on the example.  What's important is 
to pick a stratification for which a $\tsigma$ can be found that respects
the stratification, for which $H^*(S_0)$ and its limit $\check H^*(\Xi_0)$
are computable, and for which the quotient spaces $S_k/S_{k-1}$ are 
well-behaved.

\subsection{Eventual ranges}

The homotoped substitution $\tsigma$ maps $S_k$ to itself, but this
map need not be onto. Since $\tsigma$ is cellular and $S_k$ consists of a
finite number of cells, the nested sequence of spaces $S_k \supset
\tsigma(S_k) \supset \tsigma^2(S_k) \supset \cdots$ eventually stabilizes to
the {\em eventual range}, which we denote $(S_k)_{ER}$. 

The key algebraic facts about eventual ranges are:
\begin{theorem} Let $S_0 \subset S_1 \subset \cdots \subset S$ be a nested
sequence of finite CW complexes, and let $\tsigma$ be a map that sends
each $S_k$ to itself. Let $(S_k)_\ER$ denote the eventual range of $S_k$
and $\Xi_k$ the inverse limit of $S_k$. 
\begin{enumerate}
\begin{item}
$\Xi_k = \inv ((S_k)_\ER, \tsigma)$
\end{item}\begin{item}
$\check H^k(\Xi_k) = \dir (H^k((S_k)_\ER), \tsigma^*)$  
\end{item}\begin{item}
$\check H^k(\Xi_k, \Xi_{k-1})
= \dir (H^k((S_k)_\ER, (S_{k-1})_\ER), \tsigma^*).$
\end{item}\begin{item}
$\check H^k(\Xi_k, \Xi_{k-1})
= \dir (H^k((S_k)_\ER, S_{k-1} \cap (S_k)_\ER), \tsigma^*).$
\end{item}
\end{enumerate}\label{ilimthm}
\end{theorem}

\begin{proof}
Let $N_k$ be such that $\tsigma^{N_k}(S_k) = (S_k)_\ER$.
If $(x_0, x_1, \ldots)$ is a sequence of points in $S_k$, with
each $x_i = \tsigma(x_{i+1})$, then each $x_i=\tsigma^{N_k}(x_{i+N_k}) \in 
(S_k)_\ER$, so every point in the inverse limit of $S_k$ is actually
in the inverse limit of $(S_k)_\ER$. This proves the first claim, and
the second and third claims follow immediately
from the first.  The subtlety
is in the fourth, where we use the eventual range of $S_k$ but do not
use the eventual range of $S_{k-1}$. 
Suppose that $\alpha \in C^k((S_k)_\ER,
(S_{k-1})_\ER)$.  That is, $\alpha$ is a $k$-cochain that vanishes on chains
supported in $(S_{k-1})_\ER$.  However, if $c$ is a chain on $S_{k-1}$,
then $\tsigma^{N_{k-1}}(c)$ is a chain on $(S_k)_\ER$, so
$((\tsigma^*)^{N_{k-1}}\alpha)(c) = \alpha (\sigma^{N_{k-1}}(c))=0$.  
Thus the direct limit of $C^k((S_k)_\ER, (S_{k-1})_\ER)$ is the direct 
limit of $C^k((S_k)_\ER, (S_k)_\ER \cap S_{k-1})$.  The relative cohomology
of the inverse limit
is computed from the action of the coboundary on this direct limit, 
and so can be
computed either from the direct limit of $H^k((S_k)_\ER, (S_{k-1})_\ER)$
or from the direct limit of $H^k((S_k)_\ER, (S_k)_\ER \cap S_{k-1})$.  
\end{proof}

\section{Example: The chair}

\begin{figure}[ht]
\centerline{\epsfxsize=8truecm\epsfbox{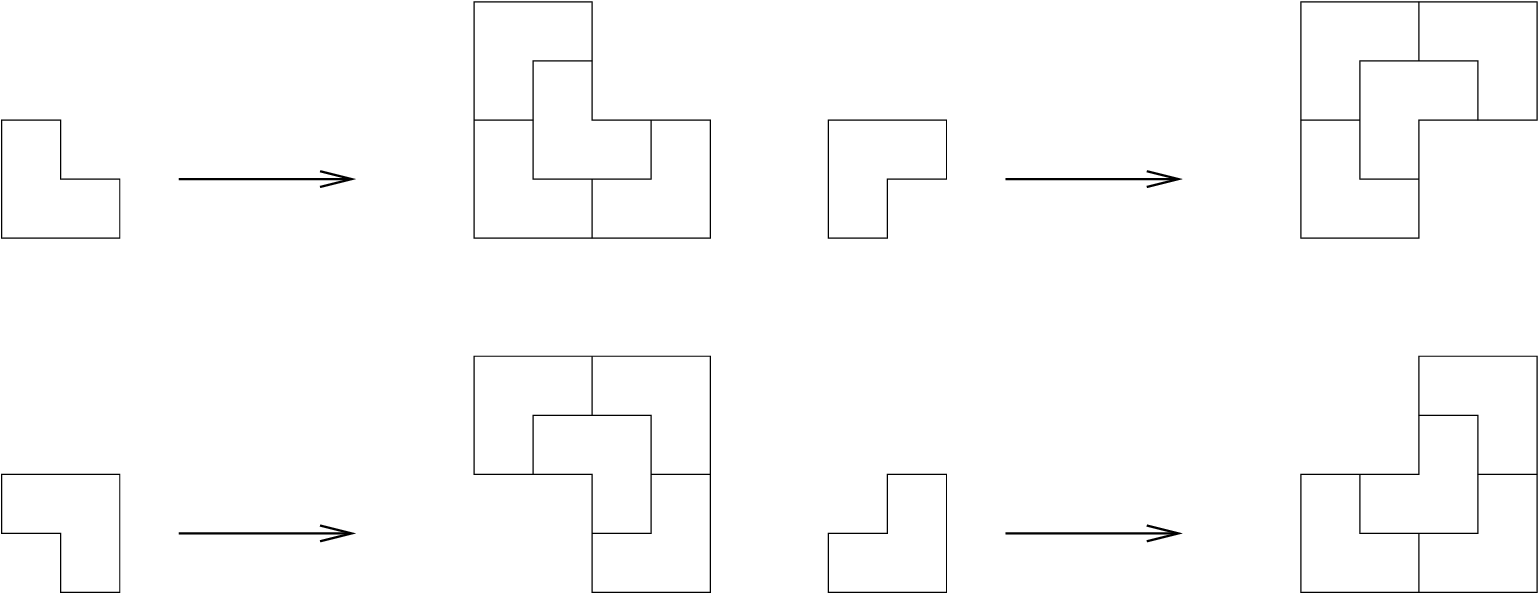}}
\caption{The chair substitution}  
\label{fig2}
\end{figure}

The chair tiling is based on a single tile, a $2 \times 2$ square with
a $1 \times 1$ corner removed, appearing in four different orientations.  
The chair tiles substitute as in Figure~\ref{fig2}. \\

If we place arrows on the three $1\times 1$ squares that make up a chair, 
$\begin{matrix} \seabox\phantom{\nwabox} \\ \neabox \nwabox \end{matrix}$,
then the chair substitution induces a substitution
$\sigma$ on arrows:

\begin{equation}\neabox \longrightarrow 
\begin{matrix} \seabox\neabox \\ \neabox\nwabox \end{matrix}, \hbox{ etc.}
\end{equation}

The conversion from chair tiles to arrow tiles and back is local, and 
the two tiling spaces are homeomorphic.  
We construct $K_t$ for the arrow substitution, and will
use this to compute the cohomology of the arrow (and therefore chair) 
tiling space.  


We use the $L^\infty$ norm on $\R^2$, for which $B_t(0)$ is geometrically
a square, of side $2t$, centered at the origin, and we pick $t<1/4$. 
We stratify $K_t$ as 
follows: $S_0$ is the set of points within $t$ of a vertex (i.e., both
the horizontal and vertical distances are less than or equal to $t$), 
$S_1$ is the set of points within $t$ of an edge, and $S_2$ is all of $K_t$.
The map $\tsigma$ is substitution 
followed by a flow towards the nearest
vertex. $S_0$ consists of one $2t \times 2t$ square for every possible 
vertex (call these ``vertex polygons''), $S_1$ consists of $S_0$ plus a
$2t \times 1-2t$ ``edge flap'' for every possible vertical edge and a
$1-2t \times 2t$ edge flap for every possible horizontal edge, and $S_2$
consists of $S_1$ plus four $1-2t \times 1-2t$ ``tile cell'' for each of 
the four possible tiles.

For typographical simplicity, we will call a northeast arrow A, a 
northwest arrow B, a southwest arrow C, and a southeast arrow D.  
Let $F_A, F_B, F_C, F_D$ 
denote the tile cells.  
There are eight vertical edge flaps, 4 of
the form $\neswabox\nwseabox$ and 4 of the form $\nwseabox\neswabox$, 
where the double-headed arrow $\neswabox$ indicates that either 
$\neabox$ or $\swabox$ can appear in this position. 
The edge flap $\neabox \nwabox$ (denoted $E_{AB}$) is glued to 
$F_A$ on the left and $F_B$ on the right, 
and hence is also glued to 
$E_{AD}$ on the left and $E_{CB}$ on the right. 
The four flaps from the configurations $\neswabox\nwseabox$
glue together to form a vertical
tube, as do the four flaps from $\nwseabox\neswabox$. 
Similarly, there are two distinct horizontal tubes capturing
allowed configurations along horizontal edges.

There are two general patterns of allowed configurations at vertices:
Type 1, of the form $\begin{matrix}\neswabox\nwseabox \\ \nwseabox\neswabox
\end{matrix}$, and
Type 2, of the form $\begin{matrix}\nwseabox\neswabox \\ \neswabox\nwseabox
\end{matrix}$.  Under substitution, each Type 1 configuration maps to the 
Type 2 configuration  $ \begin{matrix}\nwabox\neabox \\
\swabox\seabox \end{matrix}$.
Only five Type 2 corner configurations 
 are allowed for $\sigma$:   
\begin{equation} \begin{matrix}\seabox\neabox \\
\neabox\nwabox \end{matrix},
\begin{matrix}\nwabox\swabox \\
\neabox\nwabox \end{matrix},
\begin{matrix}\seabox\swabox \\
\swabox\nwabox \end{matrix},
\begin{matrix}\seabox\swabox \\
\neabox\seabox \end{matrix}, \hbox{ and }
\begin{matrix}\nwabox\neabox \\
\swabox\seabox \end{matrix}. 
\end{equation}
Each of these is taken to itself under substitution, so these
constitute the eventual range $(S_0)_{ER}$.

The diagram below illustrates a part of the complex
$K_t$.  

\begin{picture}(20, 100)(-115, -35)

\put(34,17){\line(0,-1){69}}
\put(34,17){\line(1,0){69}}
\put(34,17){\line(0,1){38}}
\put(34,17){\line(-1,0){38}}

\put(65,-14){\line(0,-1){38}}
\put(65,-14){\line(0,1){69}}
\put(65,-14){\line(1,0){38}}
\put(65,-14){\line(-1,0){69}}




\put(12, -35){$F_C$}

\put(35, -10){$\swarrow$}


\put(53, 7){$\nearrow$}
\put(80, 32){$F_A$}


\put(53, -11){$\searrow$}
\put(82, -35){$F_D$}

\put(15, 32){$F_B$}
\put(36, 7){$\nwarrow$}

\put(80,0){$E^A_D$}
\put(0,0){$E^B_C$}
\put(40,40){$E_{BA}$}
\put(40,-40){$E_{CD}$}
\dottedline{2}(49,-30)(49, 35)

\dottedline{2}(18,2)(81, 2)
\end{picture}

\bigskip

\begin{figure}
\centerline{\epsfysize=0.8truein\epsfbox{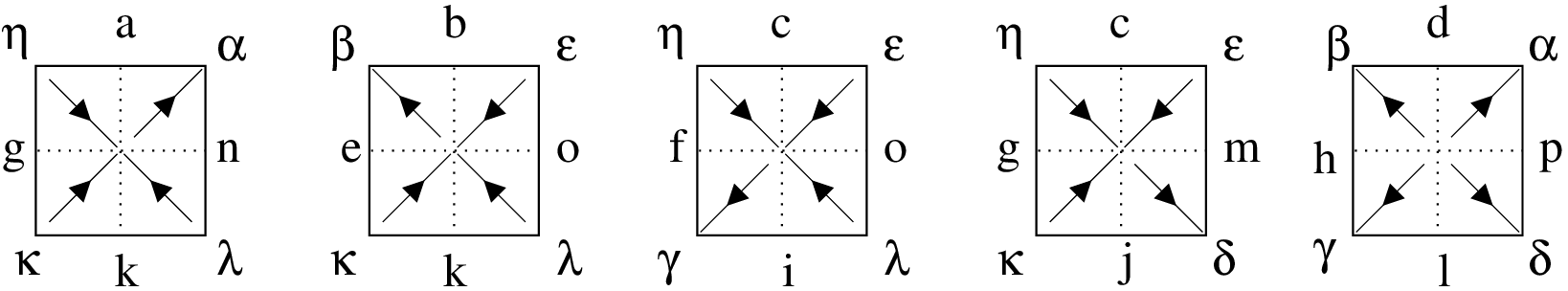}}
\caption{The 5 vertex squares in $(S_0)_\ER$}
\label{fig3}
\end{figure}

$(S_0)_\ER$ consists of five vertex squares
with boundary identifications as in Figure \ref{fig3}.  
Note that this complex has 5 faces, 16 edges and 8 vertices, for an
Euler characteristic of $-3$.  It is easy to check that 
$H^2((S_0)_\ER)=0$, and $H^0=\Z$, so $H^1 = \Z^4$. $(S_0)_\ER$ has the
homotopy type of the wedge of four circles. Since $\tsigma$ just 
permutes the cells of $(S_0)_\ER$, $\tsigma^*$ is an isomorphism and
$\check H^*(\Xi_0)=H^*((S_0)_\ER)$.  That is, 
$\check H^0(\Xi_0)=\Z$, $\check H^1
(\Xi_0)=\Z^4$ and $\check H^2(\Xi_0)=0$.

As noted earlier, the complex $S_1$ consists of $S_0$ together with 
16 edge flaps.  These edge flaps form four tubes, and $H^1(S_1, S_0)
= H^2(S_1, S_0)=\Z^4$.  
Substitution takes each edge flap of the form
$\begin{matrix} \neswabox\\ \nwseabox \end{matrix}$ 
to itself plus 
$\begin{matrix} \nwabox \\ \swabox \end{matrix}$, takes
$\begin{matrix} \nwseabox\\ \neswabox \end{matrix}$ to itself plus 
$\begin{matrix} \neabox \\ \seabox \end{matrix}$, takes
$\neswabox\nwseabox$ to itself plus $\nwabox\neabox$, and 
takes $\nwseabox\neswabox$ to itself plus $\swabox\seabox$.
$\tsigma^*$
acts by $\left ( \begin{smallmatrix} 1 & 1 & 0 & 0 \cr
1 & 1 & 0 & 0 \cr 0 & 0 & 1 & 1 \cr 0 &0 & 1 & 1 \end{smallmatrix} \right )$ 
on $H^1$, and
the direct limit is $\Z[1/2]^2$. One generator counts
horizontal edges, while the other counts vertical edges. 
However, $\tsigma^*$ acts trivially on $H^2$. 

The exact sequence of the pair $(\Xi_1, \Xi_0)$ reads:
\begin{equation} 0 \to \Z^4 \to \check H^1(\Xi_1) \to 
\Z^4 \xrightarrow{\delta} \Z^4 \to \check H^2(\Xi_1)  
\to 0. \end{equation}
The coboundary map $\delta$ is nonsingular with determinant 3, so 
$\check H^1(\Xi_1) 
=\check H^1(\Xi_1, \Xi_0)=\Z^4$ and $\check H^2(\Xi_1) = \Z_3$.

The complex $S_2 / S_1$ is just the wedge of four spheres, so 
$H^1(S_2, S_1)= 0$ and $H^2(S_2, S_1) =\Z^4$.  Under substitution,
$H^2(S_2, S_1)$ simply gets multiplied by the transpose of the
substitution matrix, and
the limit is $\Z[1/4] \oplus \Z[1/2]^2$.  (As a set, $\Z[1/4]$ is the same as
$\Z[1/2]$, but we denote it $\Z[1/4]$ to indicate that it quadruples under
substitution.)  The generator of $\Z[1/4]$ simply counts tiles, regardless
of type, while the generators of $\Z[1/2]^2$ count the vector sum of all
the arrows.  

Finally, we put it all together.  The exact sequence of
the pair $(\Xi_2, \Xi_1)$ reads 
\begin{equation} 0 \to \check H^1(\Xi_2) \to \Z[1/2]^2 \xrightarrow{\delta} 
\Z[1/4] \oplus \Z[1/2]^2 \to \check H^2(\Xi_2) \to \Z_3 \to 0.\end{equation} 
The map $\delta: \check H^1(\Xi_1) \to
\check H^2(\Xi_2,\Xi_1)$ is zero (since the boundary of 
every tile has a net of 
zero horizontal and zero vertical edges), 
so $\check H^1(\Xi_2)=\Z[1/2]^2$ and $\check H^2(\Xi_2)$ is an
extension of $\Z[1/4] \oplus \Z[1/2]^2$ by $\Z_3$.

Which extension is seen by considering the class in $\check H^2(\Omega_\sigma)$
that counts chair tiles regardless of orientation.  Since every chair tile
consists of three arrow tiles, this class is one third of the generator
of $\Z[1/4]$, so we have $\check H^2(\Omega) = \frac{1}{3} \Z[1/4] \oplus \Z[1/2]^2$
and $\check H^1(\Omega) = \Z[1/2]^2$.  

This was a long route to a simple answer.  As an Abelian group,
$\check H^2(\Omega)$ is isomorphic to the direct limit of $H^2(S_2,
S_1)$, and also to the cohomology of the inverse limit of the
uncollared Anderson-Putnam complex.  However, there is more to this
problem than the uncollared complex! The generator of the
rotationally invariant part of $\check H^2(\Omega)$ {\em cannot} be expressed
in term of uncollared arrow tiles.  The factor of 3 has to do with the
conversion from arrows to chairs, which requires information about the
neighborhood of each tile.  This collaring information is captured in
the $\Z_3$ contribution to $\check H^2(\Xi_1)$.

\section{Tilings with rotations}

We concentrate now on the case of tilings in the plane, so with $d=2$,
building into our theory the action of rotation groups. The work here
will allow us, in the final section, to compute the cohomology of the
pinwheel tiling.

\subsection{Three tiling spaces}

There are actually {\em three} tiling spaces that are associated with
a translationally finite substitution such as the chair.  The first,
denoted $\Omega^1$, is the one considered above, and is the closure of
the {\em translational} orbit of a single tiling.  For the chair
substitution, it consists of all chair tilings in which the edges are
horizontal and vertical. A finite rotation group ($\Z_4$ for the
chair, $\Z_{10}$ for the Penrose substitution) can act on $\Omega^1$,
and we can classify the terms in $\check H^*(\Omega^1)$ by how they
transform under rotation. For instance, in the chair tiling space the
$\Z[1/2]^2$ contributions to $\check H^2$ transform as a vector (i.e.,
by the matrix $\left ( \begin{smallmatrix} 0 & -1 \cr 1 & 0 \end{smallmatrix}
\right )$ for a 90
degree rotation), as do the $\Z[1/2]^2$ contributions to $\check H^1$,
while the $\frac{1}{3}\Z[1/4]$ contribution to $\check H^2$ is
rotationally invariant, as is $\check H^0 = \Z$. 

A second tiling space, denoted $\Omega^{rot}$, is the closure of the {\em
Euclidean} orbit of a single tiling. That is, $\Omega^{rot}$ is the space of all tilings
obtained by applying (orientation preserving)
rigid motions to elements of $\Omega^1$, with metric which 
stipulates that two tilings are close if one agrees with the other
in a large neighborhood of the origin, up to a small Euclidean motion.
For the chair substitution, this
gives tilings with edges pointing in arbitrary directions, not just 
vertically and horizontally. (Any one tiling will only have edges 
pointing in two perpendicular directions, but these directions are
arbitrary.) The third space, denoted $\Omega^0$, is
the set of all tilings modulo rotations about the origin.  These
are related by
\begin{equation} \Omega^{rot} / S^1 = \Omega^0 = \Omega^1/\Z_n,\end{equation}
where $n=4$ for the chair tiling and $n=10$ for the Penrose tiling. 

\subsection{Modified complexes and inverse limits}

We modify our earlier construction to write $\Omega^0$ and
$\Omega^{rot}$ as inverse limits of approximants $K^0_t$ and
$K^{rot}_t$. Let $E^{rot}$ be the trivial $\R^2$ bundle over
$\Omega^{rot}$, consisting of a copy of $\R^2$ for every tiling in
$\Omega^{rot}$. There are two obvious equivalence relations on
$E^{rot}$. The first says that $x \in T$ and $y\in T'$ are equivalent
($x \sim_t y$) if $[B_t(x)]_T = [B_t(y)]_{T'} + x-y$. That is, if
there is a translation that takes $[B_t(x)]_T$ to $[B_t(y)]_{T'}$ and
takes $x$ to $y$.  The second equivalence relation is that $x
\sim^{rot}_t y$ if there is a {\em rigid motion} that takes
$[B_t(x)]_T$ to $[B_t(y)]_{T'}$ and takes $x$ to $y$. $K^{rot}_t$ is
the quotient of $E^{rot}$ by $\sim_t$, while $K^0_t$ is the quotient
of $E^{rot}$ by $\sim^{rot}_t$.

As before, there are maps from $\Omega^{rot}$ and $\Omega^0$ to
$K^{rot}_t$ and $K^0_t$, taking a tiling to the equivalence class of
the origin. Viewed in this way, a point in $K^{rot}_t$ defines a
tiling on $B_t(0)$, while a point in $K^0_t$ defines a tiling on
$B_t(0)$ up to rotation about the origin. Substitution sends $K^{rot}_t$
to itself and $K^0_t$ to itself, and the inverse limits define tilings,
or tilings modulo rotation, on all of $\R^2$. In other words, 

\begin{theorem} $\Omega^{rot} = \inv (K^{rot}_t, \sigma)$ and 
$\Omega^{0} = \inv (K^{0}_t, \sigma)$.
\end{theorem}

If $\sigma$ has finite local complexity, then $K^{rot}_t$ is a
branched 3-manifold. Around any point is a 3-disk neighborhood
obtained by applying rigid motions to the patch $[B_t(0)]$. There are
branches whenever the closed ball $\bar B_t(0)$ intersects tiles that
the open ball $B_t(0)$ does not.  $K^0_t$ is a branched 2-orbifold,
with cone singularities at points that represent patches with
(finite!) rotational symmetry about the origin. For a generic choice
of $t$, these singularities are bounded away from the branch locus,
and we will henceforth assume that $t$ is so chosen.

Note that $K^{rot}_t$ has the structure of a (branched) Seifert manifold. 
Let $P : K^{rot}_t \to K^0_t$ be the natural projection. For any 
$x \in K^0_t$, $P^{-1}(x)$ is a circle in $K^{rot}_t$. If $x$
represents a patch without rotational symmetry about the origin,
then this circle is obtained
by rotating any representative of $x$ through $2\pi$. We call this circle a 
{\em generic fibre}. The preimage of a neighborhood $D$ of $x$ is $D
\times S^1$. If $x$ represents a patch with
$n$-fold rotational  symmetry, then $P^{-1}(x)$ is a circle obtained by
rotating any representative of $x$ through $2\pi/n$. We call this an {\em
exceptional fibre}. In this case, a neighborhood $N$ of $x$ is modeled by
the cone
$B_R(0)/\Z_n$, with 
$(r,\theta)$ (in polar coordinates) identified with $(r,\theta+2\pi/n)$. 
The preimage $P^{-1}(N)$ is then $B_R(0) \times S^1 / \Z_n$, with
$(r,\theta,\phi)
\sim (r, \theta+2\pi/n, \phi-2\pi/n)$; this is a solid torus, but the
longitude is the exceptional fibre $P^{-1}(x)$.

For tilings like the pinwheel, with tiles appearing in an infinite number
of orientations, there is no distinction between $\Omega^1$ and $\Omega^{rot}$,
since the closure of a translational orbit already contains tiles pointing
in arbitrary directions. $\Omega^{rot}$ is still the inverse limit of
branched 3-manifolds $K^{rot}_t$, $\Omega^0$ is the inverse limit of branched
2-orbifolds $K^0_t$, and we will not speak of $\Omega^1$. 

\subsection{Cohomologies of the three spaces}

For translationally finite substitution tilings, we will establish
some relations between the cohomologies of $\Omega^1$ and $\Omega^0$.
For all 2-dimensional substitutions with finite local complexity, we
will establish relations between the cohomologies of $\Omega^0$ and
$\Omega^{rot}$ and demonstrate the existence of torsion in $\check
H^2(\Omega^{rot})$ for the Penrose tiling.  In the next section we
will compute $\check H^*(\Omega^0)$ for the pinwheel tiling, and use
this to compute $\check H^*(\Omega^{rot})$.

\begin{theorem} If $\sigma$ is a translationally finite recognizable
  substitution, then $\check H^*(\Omega^0, \R)$ is isomorphic to the
  rotationally invariant part of $\check H^*(\Omega^1, \R)$.
\end{theorem}

\begin{proof}
Every real-valued cochain on $K_t$ can be written as a sum of cochains that
transform according to the irreducible representations of the rotation group
$\Z_n$ that acts on $\Omega^1$ (and therefore $K_t$). The coboundary map
is equivariant with respect to rotation, so the cohomology of $K_t$ is the
direct sum of terms, one for each irreducible representation of $\Z_n$. 
However, the rotationally invariant cochains are exactly the pullbacks of 
cochains on $K^0_t$, so $H^*(K^0_t,\R)$ is the rotationally invariant
part  of $H^*(K_t,\R)$.  Since $\sigma$ commutes with rotations, the same
observation applies to the direct limits 
$\check H^*(\Omega^0,\R) = \dir (H^*(K^0_t,\R), \sigma^*)$ and
$\check H^*(\Omega^1,\R) = \dir (H^*(K_t,\R), \sigma^*)$.
\end{proof}

This argument does not work with integer coefficients, since an
integer-valued cochain cannot necessarily be written as a sum of
irreducible components. For instance, if $n=2$, and two chains are
related by rotation by $\pi$, then the cochain $(1,1)$ is rotationally
invariant and the cochain $(1,-1)$ corresponds to the nontrivial
irreducible representation of $\Z_2$, but $(1,0)$ cannot be written as an
integer linear combination of $(1,1)$ and $(1,-1)$. We do not expect the conclusion 
of the theorem to hold with integer coefficients in general, but we know as yet of no
specific counterexamples. 

\begin{theorem}\label{rat}
  If $\sigma$ is a recognizable 
substitution with finite local complexity, then the
  real cohomology of $\Omega^{rot}$ is the same as that of
  $\Omega^0 \times S^1$.
\end{theorem}

\begin{proof}
Pick a good cover $\cU$ of $K^0_t$, such that each symmetric point lies in a 
single open set that does not touch the branch locus, and such that each
open set contains at most one symmetric point. This induces a cover
$\cV$ of $K^{rot}_t$ such that every set, and every non-empty intersection of
sets, has the topology of a circle. However, it's not always the same circle! 
Over the neighborhoods of symmetric points the circle is an exceptional
fibre,  while over all other neighborhoods, and over intersections of
neighborhoods, it is a generic fibre. 

Now consider the spectral sequence of the \v Cech-de~Rham complex for the
cover $\cV$ of $K^{rot}_t$.  That is, $E^0_{p,q}$ consists of
$q$-forms on the $p+1$-fold intersections of sets in $\cV$, $d_0$ is a
de Rham differential, $d_1$ is a \Cech differential, and so on.  Since
each nonempty intersection of sets in $\cV$ 
has the topology of a circle, we get an $E^1$
term whose 0th and 1st rows are each the \Cech complex of $\cU$ and
whose other rows are zero. The generators of the first row can be viewed
as $d\theta/2\pi$ times the generators of the zeroth row, where
$d\theta$ is the angular form on the generic fibre.  The calculations
involving $d_1$ are identical on the two rows, and the $E_2$ term is
then the \Cech cohomology of $\cU$ on the zeroth row, and again on the
first row. All that remains is to compute the differential $d_2: E_{0,1}^1
\to E_{2,0}^1$. 

To do this we start with a generator $\alpha$ 
of $E_{0,1}^1=\R$, represent it as a 
1-form on each open set, take a \Cech difference of these 1-forms, write
the result $\beta$ as the exterior derivative of a function $\gamma$ on 
the intersection of sets in $\cU$, take the \Cech differential of
$\gamma$, and view it as a class in $\check H^2(\cU, \R)$. However, 
 $K^{rot}_t$ admits a closed global angular form $d\theta$:
given any two nearby tiling patterns, we can unambiguously determine the
small angle of rotation needed to make them match, up to translation. 
Picking $\alpha$ to be this angular form, $\beta$ is identically zero,
so $d_2$ is the zero map.  

This shows that the de Rham cohomology of $K^{rot}_t$ is the same as 
the de Rham cohomology of $K^0_t \times S^1$. Taking a limit under $\sigma^*$
establishes the theorem.
\end{proof}

Note that this proof depends on the ability to find a form $d\theta/2\pi$ 
that evaluates to 1 on every generic fibre and evaluates to $1/n$ on 
exceptional fibres of order $n$.  Working with integer coefficients, that
is impossible. A cochain that evaluates to an integer on an exceptional
fibre  must evaluate to a multiple of $n$ on a generic fibre. We shall
see that this typically gives rise to  torsion in $\check H^2(\Omega^{rot})$. 

\begin{theorem}\label{Intspecseq} Let $\sigma$ be a recognizable 
substitution with finite local complexity. Suppose that $\Om^0$ contains
exactly $m$ points with $n$-fold rotational symmetry and no other symmetric
points. Then there exists a spectral sequence converging to
 $\check H^*(\Om^{rot})$ whose $E^2$ term is 
$$
\hskip -1.8in
\begin{picture}(130, 120) 
\setlength{\unitlength}{2pt}
\put(0,0){\vector(0,3){55}} 
\put(0,0){\vector(3,0){130}} 
\put(30,0){\line(0,3){40}} 
\put(80,0){\line(0,3){40}} 
\put(110,0){\line(0,3){40}} 
\put(0,20){\line(3,0){110}} 
\put(0,40){\line(3,0){110}} 
\put(8,28){$\check H^0(\Om^0)$}
\put(38,28){$\check H^1(\Om^0) \oplus \Z_n^{m-1}$}
\put(88,28){$\check H^2(\Om^0)$}
\put(8,8){$\check H^0(\Om^0)$}
\put(38,8){$\check H^1(\Om^0)$}
\put(88,8){$\check H^2(\Om^0)$}
\end{picture}
$$
Furthermore, the differential $d_2$ is zero on passing to real
coefficients.
\end{theorem}

\noindent{\em Remark.} The map $d_2$ need not be integrally zero. 
If $\check H^2(\Om^0)$ has torsion, then $d_2$ can map $E^2_{0,1}=\Z$
to a finite cyclic subgroup of $\check H^2(\Om^0)$: this is exactly what happens in the calculation of
the pinwheel tiling where $d_2$ sends the generator of $E^2_{0,1}$ to a
torsion element of order 2. 

\begin{proof}[Proof of theorem]

Pick $t$ large enough that the symmetric points of 
$K^0_t$ are in 1-1 correspondence with the symmetric points of $\Om^0$, and 
pick $t$ such that the exceptional points do not lie on the branch locus. 
As in the proof of Theorem \ref{rat}, 
pick a good cover $\cU$ of $K^0_t$, such that each symmetric point
lies in a single open set, and such that each open set contains at most one
symmetric point. 
Let $\cV$ be the preimage of $\cU$ under the projection
$P: K^{rot}_t \to K^0_t$. 
Instead of considering the {\v Cech}-de Rham complex of $\cV$, consider
the {\v Cech}-singular complex, using integer coefficients throughout. 
Since every neighborhood has the homotopy type of a circle, $E^1$
consists of two rows. The bottom row is the \Cech complex of $\cU$, and
the first  row is similar, with an infinite cyclic group for every
non-empty intersection  of sets in $\cU$.  

However, the infinite cyclic groups for the first row cannot all be
identified.  Over sets that do not contain symmetric points, the
generators evaluate to 1 on the generic fibre (call these groups
$\Z$).  Over neighborhoods of symmetric points the generators
evaluate to 1 on the exceptional fibre, and hence to $n$ on the
generic fibre (call these groups $n\Z$).  Since each symmetric point
lies in just one open set, this only affects $E^1_{0,1}$. On the bottom
row, all generators count points, and all groups are identified with $\Z$. 

The computations involving $d_1$ yield the \Cech cohomology of
$K^0_t$ on the bottom row.  For $p \ge 1$, the map $d_1: E^1_{p,1} \to
E^1_{p+1,1}$ is the same as the map $d_1: E^1_{p,0} \to E^1_{p+1,0}$. 
This implies that, for $p >1$,
$E^2_{p,1}$ identifies with $E^2_{p,0}$, and that the kernel of $d_1$ in
$E^1_{1,1}$ can be identified with the kernel of $d_1$ in $E^1_{1,0}$.
However, with our identifications, $E^1_{0,1}$ corresponds only to an
index $n^m$ subgroup of
$E^1_{0,0}$. This affects both the kernel and image of $d_1: E^1_{0,1}
\to E^1_{1,1}$. The kernel is generated by a cochain that evaluates to 
$n$ on every generic fibre and to 1 on every exceptional fibre. (In our
identification of the first and second rows, this would be all multiples of
$n$ in $E^2_{0,0}=\Z$.) The image $d_1(E^1_{0,1})$ corresponds to an index
$n^{m-1}$  subgroup of $d_1(E^1_{0,0})$. As a result, $E^2_{1,1} \cong
E^2_{1,0} \oplus
\Z_n^{m-1}$. 

Since $d_2$ was zero as a map in the {\v Cech}-de Rham double
complex, it must
be  zero as a map in the {\v Cech}-singular double complex modulo torsion.
This proves the theorem at the level of approximants. 
Finally, we note that substitution sends symmetric patterns to symmetric
patterns, and all symmetric points in $K^0_t$ correspond to symmetric points
in $\Om^0$, so substitution can only permute these points. The contributions
of the exceptional fibres therefore survive to the limit, and we obtain 
the theorem as a statement about $\Om^{rot}$ and $\Om^0$. 

\end{proof} 

\noindent{\em Remark.} The torsion appearing in
$E^2_{1,1}=E^\infty_{1,1}$ can also be understood in terms of
homology. A generic fibre is  homologous to $n$ times any exceptional
fibre, but exceptional fibres are  not homologous to each other. Rather,
the difference between any two exceptional fibres is a torsion element of
order $n$, and these differences generate a $\Z_n^{m-1}$ subgroup of 
$H_1(K^{rot}_t)$. By the universal coefficient theorem, torsion in $H_1$
gives rise to torsion in $H^2$. 

We turn to some examples.

\begin{example}
{\em The chair tiling has $\check H^1(\Om^1) = \Z[1/2]^2$ and $\check
H^2(\Om^1) =
\frac{1}{3} \Z[1/4] \oplus \Z[1/2]^2$. Direct calculation, and, equivalently here, restricting to the rotationally
invariant part, gives $\check H^1(\Om^0)=0$ and $\check
H^2(\Om^0)=\frac{1}{3}\Z[1/4]$. The space $\Om^0$ contains one point of
4-fold rotational symmetry, obtained by repeatedly substituting the
pattern where four arrows point out from the origin. The spectral sequence
that computes $\check H^*(\Om^{rot})$ has thus only four non-zero terms
by the second page: $E_{0,0}^2 =E_{0,1}^2 = \Z$, and $E_{2,0}^2 =
E_{2,1}^2=\frac{1}{3}
\Z[1/4]$. Since $E_{2,0}^2$ has no torsion,
$d_2=0$, so $E^\infty = E^2$. There are no extension problems, and we
can read off $\check H^0(\Om^{rot}) = \check H^1(\Om^{rot}) = \Z$,
with the generator of $\check H^1$ evaluating to 1 on the exceptional
fibre and to 4 on a generic fibre, and $\check
H^2(\Om^{rot}) = \check H^3(\Om^{rot}) = \frac{1}{3} \Z[1/4]$. }
\end{example}

\begin{example}
{\em 
The Penrose tiling has $ \check H^1(\Omega^1)=\Z^5$ and $ \check H^2(\Omega^1)
=\Z^8$.  The rotationally invariant part of this is 
$ \check H^1(\Omega^0)=\Z$, $ \check H^2(\Omega^0)=\Z^2$. (This can also be computed
directly.)  There are two Penrose tilings with
5-fold rotational symmetry, so our spectral sequence has 
$E^2_{1,1} = \Z \oplus \Z_5$ rather than $\Z$. As $E^\infty=E^2$, and we get 
$ \check H^0(\Omega^{rot})=\Z$, $ \check H^1(\Omega^{rot})=\Z^2$, and $ \check H^3(\Omega^{rot})=\Z^2$, 
and $\check H^2(\Omega^{rot})$ fits into the short exact sequence}
\begin{equation}
0 \to \Z^2 \to \check H^2(\Om^{rot}) \to \Z \oplus \Z_5 \to 0.
\end{equation}
\end{example}

To complete this calculation we need a way to solve the
extension problem. In fact we record the argument as a general result
as we shall have recourse to it later in the next section to complete
the pinwheel calculations: in all cases like this the extension problem
splits and there is torsion in $\check H^2(\Omega^{rot})$ of rank $m-1$. In particular, for the 
Penrose tiling, $\check H^2(\Om^{rot}) = \Z^3 \oplus \Z_5$

\begin{theorem}\label{split}
Suppose $\sigma$ and $\Omega^0$ satisfy the hypotheses of Theorem
\ref{Intspecseq}, and suppose also that
$E^\infty_{2,0}=H^2(\Omega^0)/\mbox{Im}\,d_2$ is torsion free. Then 
$$\check H^2(\Omega^{rot})\cong E^\infty_{1,1}\oplus E^\infty_{2,0}\,.$$
In particular, $\check H^2(\Omega^{rot})$ has torsion subgroup $\Z_n^{m-1}$.
\end{theorem}

\begin{proof} We consider a different decomposition of the spaces $\Om^0$
and $\Om^{rot}$, one which effectively gives us a splitting map for the
torsion subgroup in the extension problem. As in the last theorem, we
work with the approximation spaces $K_t$, the final result coming from
passing to the limit under $\sigma$. From the result of Theorem \ref{Intspecseq} we know that $\check H^2(\Omega^{rot})$ has torsion subgroup a subgroup of $\Z_n^{m-1}$; it suffices to show that it is at least $\Z_n^{m-1}$.

Let $G\subset K_t^0$ denote the union of those open sets in $\cU$ which
contain an exceptional point, thus $G$ is the disjoint union of $m$ open
discs, and let $F\subset K_t^{rot}$ be its preimage. We consider the
Mayer-Vietoris decomposition of $K_t^{rot}$ as $K_t^{rot}\setminus F$ and
$\overline F$, the closure of $F$ in $K_t^{rot}$.

Up to homotopy, $\overline F$ is a disjoint union of $m$ copies of $S^1$,
identifiable with the $m$ exceptional fibres, while the intersection of
the two subspaces is the union of $m$ copies of a 2-torus, $\cT^2$. In
cohomology, $H^1(\cT^2)=\Z^2$ and we can choose coordinates so that one
copy of $\Z$ evaluates to 1 on the generic fibre (and so to $n$ on the
exceptional one), which we shall call the {\em fibre coordinate}, while
the other coordinate represents the cohomology of the boundary circle in
the corresponding boundary disc of the open set in
$\cU$.

The relevant part of the Mayer-Vietoris sequence runs
$$\cdots\to H^1(K_t^{rot}\setminus F)\oplus\Z^m\buildrel\kappa\over\longrightarrow
\Z^{2m}\buildrel\delta\over\longrightarrow H^2(K_t^{rot})\to\cdots$$
where the $\Z^m$ represents $H^1$ of the set of exceptional fibres, and
the
$\Z^{2m}$ represents $H^1$ of the boundary tori. The map $\kappa$ on each
of these $\Z$ summands maps $1$ to $(n,0)$ in the corresponding $H^1(\cT^2)$. 
On the other hand, since there are paths in $K_t^0$ between each
exceptional point, the image of
$\kappa$ restricted to $H^1(K_t^{rot}\setminus F)$ in the $\Z^m$
corresponding to the fibre coordinates of the 2-tori is always diagonal.
Thus the cokernel of $\kappa$ contains $n$-torsion of rank $m-1$
and hence there is at least a copy of $\Z_n^{m-1}$ in $H^2(K_t^{rot})$. 
\end{proof}

With care, this argument can be extended to handle cases where
$E^\infty_{2,0}$ also contains torsion.

The results of this section did not rely on the details of our
approximants $K^{rot}_t$ and $K^0_t$. They could just as well have
been derived using Anderson-Putnam approximants and their rotational
analogs \cite{ORS}. However, the spaces $K^0_t$ provide a 
powerful tool for computing $\check H^*(\Om^0)$. We illustrate this with an
example that had heretofore resisted computation, the pinwheel tiling. 

\section{Example: The pinwheel}

\begin{figure}[ht]
\centerline{\epsfysize=1.4truein\epsfbox{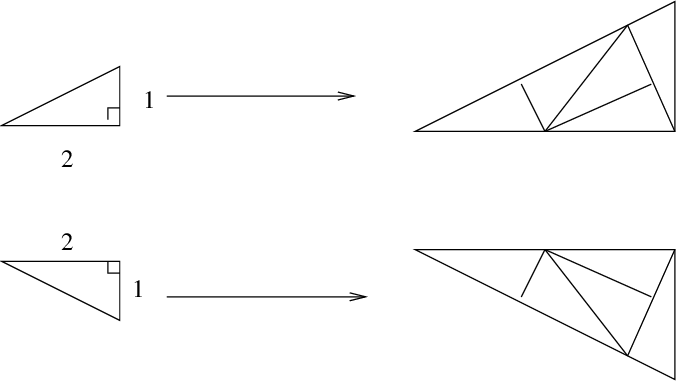}}
\caption{The pinwheel substitution}
\label{fig4}
\end{figure}

The pinwheel substitution \cite{pinwheel} involves two kinds of 
$(1,2,\sqrt{5})$ right triangles, with the substitution
rule given in Figure \ref{fig4}.

A triangle with vertices at (0,0), (2,0) and (2,1) is called right-handed,
and one with vertices at (0,0), (2,0) and (2,-1) is called left-handed.  
We call the acute-angled vertices $B_R$, $S_R$, $B_L$, and $S_L$ (for 
big-right, small-right, big-left, and small-left). We call the hypotenuses
of right-handed and left-handed tiles $H_R$ and $H_L$.  

Pick $t$ small, and stratify $K^0_t$ with $S_0$ being points within $t$
of two or more edges (i.e., close to a vertex), $S_1$ being points within
$t$ of an edge, and $S_2$ being everything. As always, the substitution
$\sigma$ does not send $S_0$ to $S_0$ or $S_1$ to $S_1$, but it is easy to
find a $\tsigma$, homotopic to $\sigma$, that does. 

Vertices in the pinwheel tiling involve combinations of big acute, small
acute, and right angles. On substitution, each right angle gets divided into
a big acute angle and a small acute angle, so configurations with right
angles do not appear in $(S_0)_\ER$.   

There are eight
kinds of vertices in $(S_0)_\ER$, corresponding to the patterns
$B_R B_L S_L S_R S_L B_L B_R S_R$,  
$B_L B_R S_R S_L S_R B_R B_L S_L$,  
$B_L B_R S_R S_L B_L B_R S_R S_L$,\\
$B_R B_L S_L S_R B_R B_L S_L S_R$,
$B_L B_R B_L B_R S_R S_L S_R S_L$, 
$B_R B_L B_R B_L S_L S_R S_L S_R$,\\
$B_R S_R S_L B_L B_R B_L S_L S_R$, and 
$B_L S_L S_R B_R B_L B_R S_R S_L$, where we list the faces counterclockwise
around the vertex.  Under substitution, the first pattern becomes the second
(and vice-versa), the third becomes the fourth, the fifth becomes the sixth,
and the seventh becomes the eighth.  We represent all but the third and fourth
as octagons, meeting the prototile faces at points and the edge flaps along
intervals.  

The third and fourth patterns have
180 degree rotational symmetry. Points related by this rotation are identified
under $\sim^{rot}_t$, so in $K^0_t$ the neighborhoods of these vertices 
correspond to quadrilaterals, rather than octagons, 
with patterns $B_L B_R S_R S_L$ and 
$B_R B_L S_L S_R$, respectively.  

There are eight edges in $(S_0)_\ER$, corresponding to the transitions
$B_R B_L$, $B_L B_R$, $S_R S_L$, $S_L S_R$, $B_R S_R$, $B_L S_L$, 
$S_R B_R$ and $S_L B_L$, and there are four vertices, namely 
$B_R$, $B_L$, $S_R$ and $S_L$.  The boundary maps are easily computed:
\begin{equation} \partial_1 = 
\begin{pmatrix} -1 & 1 & 0 & 0 & -1 & 0 & 1 & 0 \cr
1 & -1 & 0 & 0 & 0 & -1`& 0 & 1 \cr 0 & 0 & -1 & 1 & 1  & 0 & -1 & 0 \cr
0 & 0 & 1 & -1 & 0 & 1 & 0 & -1 \end{pmatrix}\end{equation}
\begin{equation} \partial_2 = \begin{pmatrix}
1 & 1 & 0 & 1 & 1 & 2 & 1 & 1 \cr
1 & 1 & 1 & 0 & 2 & 1 & 1 & 1 \cr
1 & 1 & 1 & 0 & 2 & 1 & 1 & 1 \cr
1 & 1 & 0 & 1 & 1 & 2 & 1 & 1 \cr
1 & 1 & 1 & 0 & 1 & 0 & 1 & 1 \cr
1 & 1 & 0 & 1 & 0 & 1 & 1 & 1 \cr
1 & 1 & 0 & 1 & 0 & 1 & 1 & 1 \cr
1 & 1 & 1 & 0 & 1 & 0 & 1 & 1
\end{pmatrix}
\end{equation}
and we compute $\check H^2(\Xi_0)=\check H^2((S_0)_\ER)=\Z^5$, 
$\check H^1(\Xi_0) = \check H^1((S_0)_\ER)=\Z^2$ and 
$\check H^0(\Xi_0) = \check H^0((S_0)_\ER)=\Z$.

Edges in the pinwheel tiling come in two types: hypotenuses (of length
$\sqrt{5}$) and other edges (of integer length).  Substitution interchanges
the two classes, so the contribution of each class to 
$\check H^*(\Xi_1, \Xi_0)$ must be the same, and we need only 
study the hypotenuses.

\begin{figure}[ht]
\centerline{\epsfysize=0.7truein\epsfbox{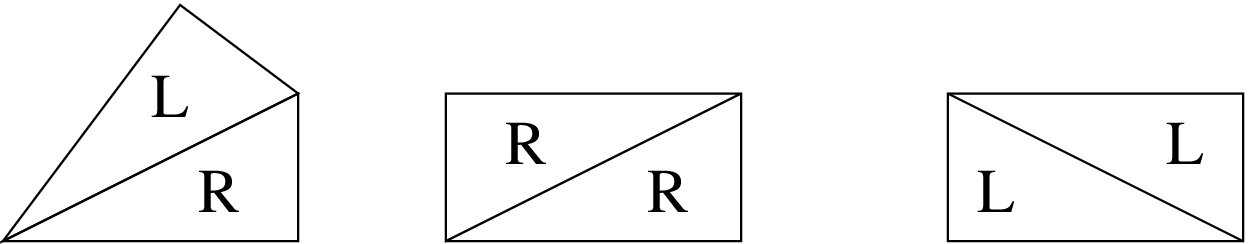}}
\caption{Three kinds of hypotenuses, one without symmetry and two with
symmetry.}
\label{fig4a}
\end{figure}

There are three kinds of hypotenuse edges, as shown in Figure~\ref{fig4a}:
those where a right-handed tile meets a left-handed tile, those where two
right-handed tiles meet, and those where two left-handed tiles meet. 
The first edge flap (call it $A$) is a $2t \times (\sqrt{5}-2t)$ rectangle
running along the hypotenuse. 
Since the other configurations have rotational symmetry, 
the second and third edge flaps (call them $B$ and $C$) 
are quotients of $2t \times (\sqrt{5}-2t)$ rectangles by rotation.  
The boundary of a $B$ edge flap
is just one hypotenuse of a right handed tile, not two.  Our
boundaries are $\partial(A)=H_R + H_L$, $\partial(B)=H_R$, and
$\partial(C)=H_L$, so $H^2=\Z$ and $H^1=0$.  Adding the contributions
of the other class of edges, we get $H^2((S_1)_\ER, S_0)=\Z^2$.
Furthermore, squared substitution multiplies these entries by 3, so 
$\check H^2(\Xi_1, \Xi_0)= \Z[1/3]^2$.  Combining the
edge flaps and vertex disks, the coboundary map from $\check H^1(\Xi_0)$ to
$\check H^2(\Xi_1,\Xi_0)$ is zero, so $\check H^2(\Xi_1)=\Z^5 + \Z[1/3]^2$, 
$\check H^1(\Xi_1) = \Z^2$ and $\check H^0(\Xi_1)=\Z$.

There are two prototiles, so $S_2/S_1$ is the wedge of two
spheres and $H^2(S_2, S_1)=\Z^2$,
and $H^1(S_2, S_1)=H^0(S_2, S_1)=0$.   Under substitution, $H^2$
transforms by the matrix 
$M=\left ( \begin{smallmatrix} 2 & 3 \cr 3 & 2 \end{smallmatrix} \right )$,
and the direct limit is $\check H^2(\Xi_2, \Xi_1)= \Z[1/5] \oplus \Z$.

In the long exact sequence of the pair $(\Xi_2, \Xi_1)$, the coboundary map
$\delta: \check H^1(\Xi_1) \to \check H^2(\Xi_2, \Xi_1)$ is not trivial.  
This is most readily seen at the level of approximants, as the image
of $\delta: H^1(S_1) \to H^2(S_2, S_1) = \Z^2$ is 
all multiples of $(2,-2)$, a set that is invariant under substitution. 
The cokernel is $\Z \oplus \Z_2$, and the direct limit of the 
cokernel is $\Z[1/5] \oplus \Z_2$. 
This implies that $\check H^1(\Omega^0)=\Z$ and
$\check H^2(\Omega^0)=\Z[1/5] \oplus \Z[1/3]^2 \oplus \Z^5 \oplus \Z_2$. 

There are six pinwheel tiling configurations with 180 degree
rotational symmetry.  Two correspond to the third and fourth vertex
disks.  Two correspond to the $B$ and $C$ hypotenuse edge flaps.  Two
are integer-length edge flaps obtained from the $B$ and $C$ edge flaps by
substitution.  These six configurations yield six exceptional fibres of
the fibration $\Omega^{rot} \to \Omega^0$. 

This means that the spectral sequence computing $\check H^*(\Omega^{rot})$
has $E^2$ term 
$$
\hskip -2in
\begin{picture}(130, 120) 
\setlength{\unitlength}{2pt}
\put(0,0){\vector(0,3){55}} 
\put(0,0){\vector(3,0){130}} 
\put(20,0){\line(0,3){40}} 
\put(45,0){\line(0,3){40}} 
\put(120,0){\line(0,3){40}} 
\put(0,20){\line(3,0){120}} 
\put(0,40){\line(3,0){120}} 
\put(8,28){$\Z$}
\put(23,28){$\Z \oplus \Z_2^5$}
\put(48,28){$\Z[1/5]\oplus\Z[1/3]^2\oplus \Z^5 \oplus \Z_2$}
\put(8,8){$\Z$}
\put(23,8){$\Z$}
\put(48,8){$\Z[1/5]\oplus\Z[1/3]^2\oplus \Z^5 \oplus \Z_2$}
\end{picture}
$$
The generator of $E^2_{0,1}$ is a cochain that evaluates to 1 
on every exceptional fibre and to 2 on every generic fibre.

The next step is computing $d_2$. Since $d_2$ is zero over $\R$, it must
send the generator of $E^2_{0,1}$ either to zero or to the unique torsion
element in $E^2_{2,0}$. In the first instance, $E^\infty_{0,1}$ would be
generated by a cochain that evaluates to 1 on each exceptional fibre. 
In the second instance, $E^\infty_{0,1}$ would be generated by a cochain
that evaluates to 2 on each exceptional fibre. We claim that the first
is impossible, and that the second is correct. 

To see this, we construct a path $\gamma$ in $K^{rot}_t$ such that
$2\gamma$ is homologous to an exceptional fibre. This implies that 
any cohomology class in $H^1(K^{rot}_t)$, evaluated on the exceptional
fibre, must yield an even result, and in particular cannot yield 1. Since
every cohomology class on $\Om^{rot}$ can be represented by a class in 
an approximant, there are no classes in $\check H^1(\Om^{rot})$ (or in 
$E^\infty_{0,1}$) that
evaluate to 1 on the exceptional fibre. 

\begin{figure}[ht]
\vskip -1.2truein
\centerline{\hskip 1truein \epsfysize=4truein\epsfbox{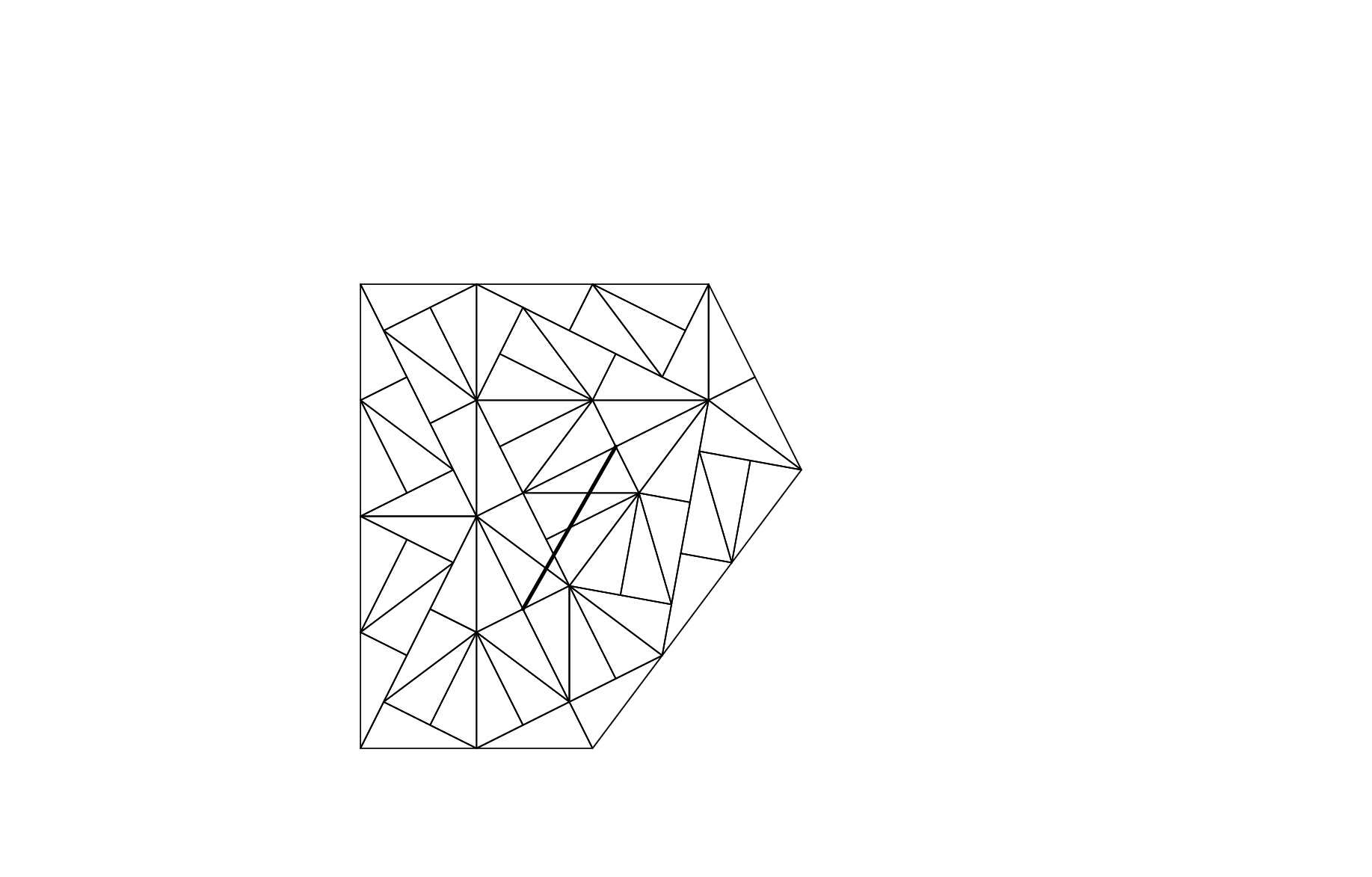}}
\vskip -0.6truein
\caption{A path in a pinwheel tiling}
\label{pinpatch}
\end{figure}

Consider the patch of the pinwheel tiling shown in Figure
\ref{pinpatch}, and the path shown in it. This path induces a closed
loop $\gamma_0$ in $K^0_t$.  Although $\gamma_0$ is not a boundary,
$2\gamma_0$ is. The loop $\gamma$ in $K^{rot}_t$ is induced by the
path in the figure, followed by a 90 degree counterclockwise rotation.
$2\gamma$ is homotopic to parallel transport by $2\gamma_0$ followed
by a 180 degree rotation about the endpoint, which in turn in
homotopic to just a 180 degree rotation about the endpoint, which is
one lap around an exceptional fibre.

$E^\infty$ is thus
$$
\hskip -2in
\begin{picture}(130, 130) 
\setlength{\unitlength}{2pt}
\put(0,0){\vector(0,3){55}} 
\put(0,0){\vector(3,0){130}} 
\put(20,0){\line(0,3){40}} 
\put(45,0){\line(0,3){40}} 
\put(120,0){\line(0,3){40}} 
\put(0,20){\line(3,0){120}} 
\put(0,40){\line(3,0){120}} 
\put(8,28){$\Z$}
\put(23,28){$\Z \oplus \Z_2^5$}
\put(48,28){$\Z[1/5]\oplus\Z[1/3]^2\oplus \Z^5 \oplus \Z_2$}
\put(8,8){$\Z$}
\put(23,8){$\Z$}
\put(48,8){$\Z[1/5]\oplus\Z[1/3]^2\oplus \Z^5$}
\end{picture}
$$
The generator of $E^\infty_{0,1}$ is a cochain that evaluates to 2 
on every exceptional fibre and to 4 on every generic fibre. 

In general, $E^\infty$ does not uniquely determe $\check H^*$.
Rather, $\check H^k(\Omega^{rot})$ fits into the exact sequence
\begin{equation} 0 \to E^\infty_{k,0} \to \check H^k(\Omega^{rot}) 
\to E^\infty_{k-1,1} \to 0.
\end{equation}
Since $E^\infty_{0,1}$ is free, and since $E^\infty_{3,0}$ vanishes, there
are no extension problems in computing $\check H^1(\Omega^{rot})=\Z^2$
or $\check H^3(\Omega^{rot})= \Z[1/5] \oplus \Z[1/3]^2 \oplus \Z^5 \oplus 
\Z_2$. The exact sequence involving $\check H^2(\Omega^{rot})$ reads
\begin{equation} 0 \to \Z[1/5] \oplus \Z[1/3]^2 \oplus \Z^5 \to 
\check H^2(\Omega^{rot}) \to \Z \oplus \Z_2^5\to 0,
\end{equation}
so by Theorem \ref{split}, $\check H^2(\Omega^{rot}) = \Z[1/5] \oplus
\Z[1/3]^2 \oplus\Z^6 \oplus \Z_2^5$.

This calculation helps resolve a longstanding question about variants
of the pinwheel tiling.  One can build a pinwheel variant with
an $(m,n,\sqrt{m^2+n^2})$ right triangle, where $m$ and $n$ are arbitrary
integers.  The linear expansion factor is $\sqrt{m^2+n^2}$, and it was
wondered \cite{ORS} whether pinwheel spaces with the same stretching 
factor (e.g., the $(8,1)$ and $(7,4)$-pinwheels) had homeomorphic 
tiling spaces.  

They do not.  The strata $S_0$ for the $(m,n)$ pinwheel spaces are more
complicated than for the ordinary (2,1) pinwheel space, but $S_2/S_1$ is
exactly as before.  The edge flaps give a $\Z^2$ contribution to
$H^2(S_2,S_1)$, and these terms scale by $m^2+n^2 -2|m-n|$ under squared
substitution.  This gives  a contribution of $\Z[1/51]^2$ to $\check H^2$
of the $(8,1)$-pinwheel and 
$\Z[1/59]^2$ to $\check H^2$ of the $(7,4)$-pinwheel.  
There is also a contribution
of $\Z[1/(m^2+n^2-2(m+n))]$ from the direct limit of the substitution matrix.
This gives a $\Z[1/47]$ contribution to $\check H^2$ of the $(8,1)$-pinwheel and
$\Z[1/43]$ to $\check H^2$ of the $(7,4)$-pinwheel.  
$H^*((S_0)_\ER)$ is invariant under substitution, so nothing in 
$\check H^*(\Xi_0)$ can cancel or mimic these $p$-adic terms.

\bigskip\noindent{\bf Acknowledgments.} All four authors thank the
Banff International Research Station and the participants in tilings
workshops held there in 2006 and 2008. Many of the ideas for this
article were developed at these workshops. 
J.H. thanks the Royal Society for support and the University
of Leicester for study leave. The work of L.S. is partially supported by 
the National Science Foundation.

\bigskip

{\Small {\parindent=0pt Department of Mathematics, Montana State University,
Bozeman, MT 59717, USA \\
barge@math.montana.edu \\
\phantom{asdf} \\
Department of Mathematics, College of Charleston, Charleston, SC 29424,
USA
\\ diamondb@cofc.edu \\
\phantom{asdf} \\
Department of Mathematics, University of Leiceister, Leicester, LE1 7RH,
England \\ j.hunton@mcs.le.ac.uk \\
\phantom{asdf} \\
Department of Mathematics, University of Texas, Austin, TX 78712, USA \\
sadun@math.utexas.edu
}}

\end{document}